\documentclass[11pt]{article}
\usepackage{anysize}
\marginsize{2.0cm}{2.0cm}{2.0 cm}{2.0 cm}
\usepackage{setspace}
\usepackage{graphicx}
\usepackage{color}
\usepackage{epstopdf}
\usepackage[titletoc]{appendix}
\usepackage{graphicx}
\usepackage{latexsym,amsfonts,amsmath,amssymb}
\newtheorem{remark}{Remark}

\newtheorem{theorem}{Theorem}
\newtheorem{lemma}{Lemma}

\newtheorem{definition}{Definition}
\newtheorem{proposition}{Proposition}
\newtheorem{example}{Example}

\def\R{\mathbb{R}}
\def\I{\mathbb{I}}
\def\PP{\mathsf{P}}
\def\EE{\mathsf{E}}

\begin{document}
\title{On the continuity of  the projection mapping from strategic measures to occupation measures in absorbing Markov decision processes}

\author{Alexey Piunovskiy\thanks{Department of Mathematical Sciences, University of
Liverpool, Liverpool, U.K.. E-mail: piunov@liv.ac.uk.}~ and Yi
Zhang \thanks{Corresponding author. School of Mathematics, University of
Birmingham, Birmingham, B15 2TT, U.K.. E-mail: y.zhang.29@bham.ac.uk.}}
\date{}

\maketitle

\par\noindent{\bf Abstract:} In this paper, we prove the following assertion for an absorbing Markov decision process (MDP) with the given initial distribution, which is also assumed to be semi-continuous: the continuity of the projection mapping from the space of strategic measures to the space of occupation measures, both endowed with their weak topologies, is equivalent to the MDP model being uniformly absorbing. An example demonstrates, among other interesting scenarios, that for an absorbing (but not uniformly absorbing) semi-continuous MDP with the given initial distribution, the space of occupation measures can fail to be compact in the weak topology.
\bigskip

\par\noindent {\bf Keywords:}  Markov decision processes. Absorbing model. Continuity. Projection mapping.
\bigskip

\par\noindent
{\bf AMS 2020 subject classification:} 90C40, 60J05

\section{Introduction}\label{s1_0}
In this paper, we consider a Markov decision process (MDP) with a Borel state space ${\bf X}$ and a Borel action space ${\bf A}$, both being endowed with their Borel $\sigma$-algebras. If there is an isolated absorbing state say $0$ in the state space, then the MDP model is called absorbing for a given initial distribution ${\PP_0}$ if under each strategy, the expected hitting time to $0$ is finite.

In terms of occupation measures, understood as the total state-action frequencies on $({\bf X}\setminus \{0\})\times {\bf A}$, endowed with the product $\sigma$-algebra, an MDP is absorbing for the given initial distribution if the occupation measure of each strategy is a finite measure. Occupation measures are important for the study of optimal control problems of MDPs with total cost criteria because the performance measure can be written as an integral of the cost function with respect to occupation measures. This turns the original MDP problem to a static optimization problem in the space ${\cal D}$ of occupation measures. Since for absorbing MDP models, ${\cal D}$ contains only finite measures, it is natural, as we do in this paper, to endow it with the weak topology, which is metrizable.

The wonderful and insightful paper \cite{b114}, with many new ideas, intended to develop a rich theory for absorbing MDPs, in particular, for the occupation measures in such MDP models. A key property was used there, which asserts that if the MDP model is absorbing and semi-continuous, then the space of occupation measures is compact with respect to the weak topology, see \cite[Lemma 4.7]{b114}. Here an MDP model is called semi-continuous if the action space ${\bf A}$ is compact, and the transition kernel $p(dy|x,a)$ is either set-wise continuous in $a\in {\bf A}$ for each $x\in {\bf X}$, or it is continuous with respect to the weak topology in $(x,a)\in{\bf X}\times {\bf A}$.

In order to prove this result, \cite{b114} exploited the following fact: for semi-continuous MDP models, the space of strategic measures is compact in weak topology, as established in \cite{b18}, and the projection mapping $O$, carrying the strategic measure of a strategy to the occupation measure of the same strategy, is continuous.

In the present paper, by means of an example, see Example \ref{PiunovskiyZhang2023cExample01} below, we demonstrate that the aforementioned assertion in \cite[Lemma 4.7]{b114} is inaccurate. This is due to the fact that in general, the projection mapping $O$ may be not continuous for absorbing semi-continuous MDP models, see Theorem \ref{PiunovskiyZhang2023cThm01}.

Our second contribution is as follows. We show that for an absorbing MDP model, provided that it is semi-continuous, the projection mapping $O$ is continuous if and only if when the MDP model is uniformly absorbing. The latter requires that, for the given initial distribution, the expected hitting time to $0$ converges uniformly with respect to all strategies. In fact, the sufficiency part follows from the same reasoning as in the proof of \cite[Lemma 4.7]{b114}, by avoiding the minor error therein. This characterization of the continuity of the projection mapping $O$ shows that the gap in \cite{b114} can be naturally closed if one further requires the MDP model there to be uniformly absorbing.

To the best of our knowledge, this definition of uniformly absorbing MDPs for a given initial distribution seems to appear in \cite{b16}. In that paper, which focused on non-atomic MDP models, its equivalence to the continuity of the projection mapping $O$ was not discussed.

A trivial example of uniformly absorbing MDP models is given by those discounted MDP models with a nonnegative discount factor strictly smaller than $1$.

If the initial distribution is concentrated on a single state, and the state space is countable or finite, endowed with the discrete topology, then the MDP model, assumed to be semi-continuous, is uniformly absorbing if it possesses a uniform Lyapunov function, see Definition \ref{d2_01} below. MDP models with a uniform Lyapunov function were studied intensively in e.g., \cite{b13,Cavazos:1992,Hordijk:1975}. In particular, for such MDP models, for every initial state, the MDP model is uniformly absorbing, see Proposition \ref{t2_04}. However, if one considers a fixed initial distribution, then the existence of a uniform Lyapunov function does not imply the model to be uniformly absorbing or even absorbing. This is demonstrated in Example \ref{PiunovskiyZhang2023cExample02}, and is sometimes overlooked. Similar models for MDPs with Borel state and action spaces were considered in e.g., \cite{Collins:2006} and \cite[Chapter 9]{b12}, which assumed that the total value of occupation measures is bounded or $w$-bounded over all initial states and strategies, and focused on the optimality equation.

Let us mention that, apart from providing an absorbing semi-continuous MDP model with the space of occupation measures being not compact, Example \ref{PiunovskiyZhang2023cExample01} incidentally demonstrates several other scenarios of interest. E.g., it also demonstrates that an absorbing MDP model for a fixed initial state may not possess a uniform Lyapunov function, and that some known solvability conditions are important.

The rest of this paper is organized as follows. We describe the MDP model in Section \ref{PiunovskiyZhang2023cSectionModel}, and recall some known facts in Section \ref{PiunovskiyZhang2023cSectionFacts}. The main results are presented in Section \ref{PiunovskiyZhang2023cSectionMain}. The main example, i.e., Example \ref{PiunovskiyZhang2023cExample01}, is formulated in Section \ref{PiunovskiyZhang2023cSectionExample}, where we also discuss it in the context of optimal control problem of MDPs. The proofs of all the statements (except those known ones) are in Section \ref{PiunovskiyZhang2023cSectionProof}. Finally, this paper is finished with a conclusion in Section \ref{PiunovskiyZhang2023cSectionConclusion}.

\section{Model description}\label{PiunovskiyZhang2023cSectionModel}

Before describing the MDP model, let us fix some notation and conventions used throughout this paper. If a space $\bf Y$ is discrete (with the discrete topology) and $\mu$ is a measure on ${\cal B}({\bf Y})$, then, for singletons, we use notation $\mu(y)$, not $\mu(\{y\})$. Integral of a function $c(\cdot)$ with respect to a measure $\mu$ is written as $\int_{\bf Y} c(y)d\mu(y)$ or $\int_{\bf Y} c(y)\mu(dy)$. In $\R=(-\infty,\infty)$, unless stated otherwise, the usual Euclidean topology is fixed.   $\delta_a(dx)$ is the Dirac measure concentrated at the point $a$, provided that the singleton $\{a\}$ is measurable, and $\I\{\cdot\}$ is the indicator function.

We fix the following primitives of a Markov decision process (MDP) model.
\begin{itemize}
\item $\bf X$ and $\bf A$ are the state and action spaces, both assumed to be nonempty topological Borel spaces. Here a topological Borel space is a Borel subset of some completely metrizable separable space. We endow ${\bf X}$ and ${\bf A}$ with their Borel $\sigma$-algebras ${\cal B}({\bf X})$ and ${\cal B}({\bf A})$, respectively.
\item The transition probability $p(dy|x,a)$ is a (measurable) stochastic kernel on ${\cal B}({\bf X})$ given ${\bf X}\times{\bf A}$.
\end{itemize}
Given the collection $\{{\bf X},{\bf A}, p\}$, given each strategy (see Definition \ref{PiunovskiyZhang2023cDef05} below) and initial distribution on $\bf X$, one can construct a probability space, and define the controlled (state) process $\{X_t\}_{t=0}^\infty$ and controlling (action) process $\{A_t\}_{t=1}^\infty$ thereon. Derman in his classic book \cite[p.4]{Derman:1970} termed the bivariate process $\{(X_t,A_{t+1})\}_{t=0}^\infty$ a Markov decision process. This is why we refer to $\{{\bf X},{\bf A},p\}$ as the MDP model.  The aforementioned construction is described next.

The space of trajectories (or say histories) is the following countable product
  \begin{eqnarray*}
  {\bf H}:={\bf X}\times({\bf A}\times{\bf X})^\infty.
  \end{eqnarray*}
The generic notation for an element of ${\bf H}$ is $\omega=(x_0,a_1,x_1,\ldots)\in {\bf H}$. We endow ${\bf H}$ with the product $\sigma$-algebra, which is also the Borel $\sigma$-algebra on it.
All the random variables, like $X_t$ and $A_{t+1}$, are just measurable mappings defined on $\bf H$:
\begin{eqnarray*}
X_t(\omega)=x_t, A_{t+1}(\omega)=a_{t+1},~\forall~t\ge 0.
\end{eqnarray*}

\begin{definition}\label{PiunovskiyZhang2023cDef05}
A control strategy, or simply say a strategy, $\pi=\{\pi_t\}_{t=1}^{\infty}$ is a sequence of stochastic kernels $\pi_t(da|x_0,a_1,x_1,\dots, a_{t-1},x_{t-1})$ on ${\cal B}({\bf A})$ given ${\bf H}_{t-1}:=({\bf X}\times{\bf A})^{t-1}\times{\bf X}$.
A strategy $\pi$ is called Markov, if for each $t\ge 1$, there is a stochastic kernel $\pi^M_t$ on ${\cal B}({\bf A})$ given ${\bf X}$ such that $\pi_t(da|x_0,a_1,x_1,\ldots, a_{t-1},x_{t-1})=\pi^M_{t}(da|x_{t-1})$.
If, for some stochastic kernel $\pi^s$ on ${\cal B}({\bf A})$ given ${\bf X}$, $\pi_t(da|x_0,a_1,x_1,\ldots, a_{t-1},x_{t-1})=\pi^s(da|x_{t-1})$ for any $t=1,2,\dots$,  then the strategy $\pi=\{\pi_t\}_{t=1}^\infty$ is called stationary, and is identified with and denoted as $\pi^s$. If a stationary strategy  $\pi^s$ takes the form $\pi^s(da|x)=\delta_{\varphi(x)}(da)$ for some measurable mapping $\varphi$ from ${\bf X}$ to ${\bf A}$, where $\delta_{\varphi(x)}(da)$ denotes the Dirac measure concentrated on the singleton $\{\varphi(x)\}$, then the strategy $\pi^s$ is called deterministic stationary, and is identified with and denoted by the underlying measurable mapping $\varphi$.
\end{definition}
The set of all strategies is denoted as $\Delta^{\rm All}$, and the set of Markov strategies is denoted as $\Delta^{\rm Markov}$.

Given a strategy $\pi=\{\pi_t\}_{t=1}^{\infty}$, for a fixed state-action pair $(z,b)\in {\bf X}\times {\bf A}$, we define its shifted strategy ${}^{(z,b)}\pi=\{{}^{(z,b)}\pi_t\}_{t=1}^\infty$ by
\begin{eqnarray}\label{PiunovskiyZhang2023cEqn21}
{}^{(z,b)}\pi_t(da|y_0,a_1,\dots,a_{t-1},y_{t-1}):=\pi_{t+1}(da|z,b,y_0,a_1,\dots,a_{t-1},y_{t-1})
\end{eqnarray}
for $y_i\in {\bf X}$  with $i\in\{0,1,\dots,t-1\}$ and $a_i\in {\bf A}$ with $i\in\{1,2,\dots,t-1\}$.

Let the initial distribution $\PP_0$ on $({\bf X},{\cal B}({\bf X}))$ be given. If a strategy $\pi$ is also fixed, then the strategic measure on $({\bf H},{\cal B}(\bf H))$, constructed in the standard way using the Ionescu-Tulcea Theorem, is denoted as $\PP^\pi_{\PP_0}$. It is the unique probability measure on $({\bf H},{\cal B}({\bf H}))$ such that
\begin{eqnarray*}
&&\PP^\pi_{\PP_0}(X_0\in dx)=\PP_0(dx);\\
&& \PP^\pi_{\PP_0}(A_{t+1}\in da|X_0,A_1,\dots,X_t)=\pi_{t+1}(da|X_0,A_1,\dots,X_t);\\
&& \PP^\pi_{\PP_0}(X_{t+1}\in dx|X_0,A_1,\dots,X_t,A_{t+1})=p(dx|X_t,A_{t+1})~\forall~t\ge 0.
\end{eqnarray*}

The corresponding mathematical expectation is denoted as $\EE_{\PP_0}^\pi$. If the initial distribution $\PP_0(dx)=\delta_{x_0}(dx)$ is degenerate, we use notations $\PP^\pi_{x_0}$ and $\EE^\pi_{x_0}$.

We denote by
\begin{eqnarray*}
{\cal P}:=\{\PP_{\PP_0}^\pi:\pi\in \Delta^{\rm All}\}
\end{eqnarray*}
the space of all strategic measures with the given initial distribution $\PP_0$.

An important class of MDP models is the absorbing model defined as follows.
\begin{definition} \label{d3}
The MDP is called absorbing (at $0$) for the given initial distribution $\PP_0$ if there is an isolated state, say $0$ in ${\bf X}$ such that $0$ is absorbing, i.e.,
\begin{eqnarray*}
p(\{0\}|0,a)\equiv 1,
\end{eqnarray*}
and
\begin{eqnarray*}
\EE^\pi_{\PP_0}[T_0]<\infty,~\forall~\pi\in\Delta^{\rm All},
\end{eqnarray*}
where
\begin{eqnarray*}
T_0:=\inf\{t\ge 0:~X_t=0\}
\end{eqnarray*}
denotes the hitting time to the state $0$ by the controlled process. As usual, $\inf\emptyset:=\infty.$
\end{definition}
Verbally, an MDP model is absorbing (at $0$) for the initial distribution $\PP_0$ if under each strategy, the expected hitting time to the isolated absorbing state $0$ is finite. The above definition of an absorbing MDP model for the initial distribution $\PP_0$ is the same as the one given by Feinberg and Rothblum in \cite[p.7]{b114}.  Given that the state $0$ is absorbing, this definition also coincides with the one in \cite[Definition 7.1]{b13}. In general, in \cite[Definition 7.1]{b13}, for the model to be absorbing for $\PP_0$, the state $0$ was not required to be absorbing itself, and it was required that $\EE^\pi_{\PP_0}[\tau_0]<\infty,~\forall~\pi\in\Delta^{\rm All},$ where
\begin{eqnarray*}\tau_0:=\inf\{t\ge 1:~X_t=0\}
\end{eqnarray*}
is the return time to state $0$.

Thus, for an absorbing MDP model (at state $0$) for the initial distribution $\PP_0$, for each strategy, the series $\sum_{t=0}^{\infty} \EE_{\PP_0}^\pi[I\{T_0>t\}]$ converges, because $\sum_{t=0}^{\infty} \EE_{\PP_0}^\pi[\I\{T_0>t\}]=\EE_{\PP_0}^\pi[T_0]<\infty$. If we require the convergence of the above series to be uniform with respect to all strategies $\pi\in\Delta^{\rm All}$, then the resulting model will be called uniformly absorbing, formulated in the next definition.
\begin{definition}\label{PiunovskiyZhang2023cDef01}
An absorbing (at 0) MDP model for the initial distribution $\PP_0$ is called {\sl uniformly absorbing} for $\PP_0$  if
\begin{eqnarray*}
\lim_{n\to\infty}\sup_{\pi\in\Delta^{\rm All}} \EE^\pi_{\PP_0}\left[\sum_{t=n}^\infty \I\{t<T_0\}\right]=0.
\end{eqnarray*}
\end{definition}
This definition was given in \cite[Definition 3.6]{b16}. Since the state $0$ was required to be absorbing in Definition \ref{d3}, the requirement in Definition \ref{PiunovskiyZhang2023cDef01} is the same as
\begin{eqnarray}\label{PiunovskiyZhang2023cEqn05}
\lim_{n\to\infty}\sup_{\pi\in\Delta^{\rm All}} \EE^\pi_{\PP_0}\left[\sum_{t=n}^\infty \I\{X_t\in({\bf X}\setminus\{0\})\}\right]=0
\end{eqnarray}
because
\begin{eqnarray*}\EE^\pi_{\PP_0}\left[\sum_{t=n}^\infty \I\{t<T_0\}\right]= \EE^\pi_{\PP_0}\left[\sum_{t=n}^\infty \I\{X_t\in({\bf X}\setminus\{0\})\}\right],~\forall~n\ge 0.
\end{eqnarray*}

This equality, when specialized to $n=0,$ allows us to formulate the definition of absorbing MDP models in terms of the finiteness of the occupation measures, defined as follows.

 \begin{definition}\label{d2_1}
Consider an MDP model $\{{\bf X},{\bf A},p\}$ with $0\in {\bf X}$ being an isolated absorbing state, i.e., $p(\{0\}|0,a)\equiv 1.$
The occupation measure $\eta^\pi_{\PP_0}$ of a strategy $\pi$ for the given initial distribution $\PP_0$ is the $[0,\infty]$-valued measure on $(({\bf X}\setminus\{0\})\times {\bf A},{\cal B}(({\bf X}\setminus\{0\})\times{\bf A}))$ defined by the formula
  \begin{eqnarray*}
  \eta^\pi_{\PP_0}(\Gamma^{\bf X}\times\Gamma^{\bf A}):=\sum_{t=1}^\infty \PP_{\PP_0}^\pi\{X_{t-1}\in\Gamma^{\bf X}, A_{t}\in\Gamma^{\bf A}\},~\forall~\Gamma^{\bf X}\in{\cal B}({\bf X}\setminus \{0\}),~\Gamma^{\bf A}\in{\cal B}({\bf A}).
  \end{eqnarray*}
\end{definition}
If the initial distribution is a singleton, say $x_0,$ then we write $\eta^\pi_{x_0}$ for the occupation measure of a strategy $\pi.$
Let
\begin{eqnarray*}
{\cal D}:=\{\eta^\pi_{\PP_0}:\pi\in\Delta^{\rm All}\}
\end{eqnarray*}
be the space of all occupation measures for the initial distribution $\PP_0$.

Now we see that an MDP model with an absorbing state $0$ is absorbing (at $0$) for the initial distribution $\PP_0$ if and only if
\begin{eqnarray*}
\eta^\pi_{\PP_0}(({\bf X}\setminus \{0\})\times {\bf A})<\infty
\end{eqnarray*}
for each $\pi\in\Delta^{\rm All}.$

For brevity, we make the next definition.
\begin{definition}\label{PiunovskiyZhang2023cDef02}
An MDP model $\{{\bf X},{\bf A},p\}$ is called semi-continuous (S)  if the following two conditions are satisfied:
\begin{itemize}
\item[(a)] The action space $\bf A$ is a compact topological Borel space.
\item[(b)] For each $x\in{\bf X}$, the function $\int_{\bf X} u(y)p(dy|x,a)$ is continuous in $a\in {\bf A}$ for every bounded measurable function $u(\cdot)$.
\end{itemize}
An MDP model $\{{\bf X},{\bf A},p\}$ is called semi-continuous (W)  if the following two conditions are satisfied:
\begin{itemize}
\item[(c)] The action space $\bf A$ is a compact topological Borel space.
\item[(d)] The function $\int_{\bf X} u(y)p(dy|x,a)$ is continuous in $(x,a)\in {\bf X}\times {\bf A}$ for every bounded continuous function $u(\cdot)$.
\end{itemize}
An MDP model $\{{\bf X},{\bf A},p\}$ is called semi-continuous if it is either semi-continuous (S) or semi-continuous (W).
\end{definition}
Conditions (a,b) in Definition \ref{PiunovskiyZhang2023cDef02} are respectively called Condition (S)(1,2) in \cite[Section 6]{b18}, whereas conditions (c,d) in Definition \ref{PiunovskiyZhang2023cDef02} are respectively called Condition (W)(1,2) in \cite[Section 5]{b18}. This is why we termed semi-continuous (S) and semi-continuous (W) models in the above. If ${\bf X}$ is countable or finite, and is endowed with the discrete topology, then semi-continuity (S) and semi-continuity (W) mean the same.
Admittedly, the term of ``semi-continuous MDP models'' has other meanings in the literature. Nevertheless, in the present paper, its meaning is unambiguous.

\section{Some relevant facts}\label{PiunovskiyZhang2023cSectionFacts}

In this section we present some facts about MDP models and optimal control problems of MDPs.
\subsection{Facts about MDP model}

\subsubsection{MDP with a uniform Lyapunov function}
An important class of semi-continuous MDP models with a countable or finite state space that are uniformly absorbing (at $0$) for a given initial state $x_0$ is given by those that admit a uniform Lyapunov function, defined as follows.

\begin{definition}\label{d2_01}
Consider a semi-continuous MDP model with a countable or finite state space $\bf X$ (endowed with the discrete topology) with the isolated absorbing state $0$. A $[1,\infty)$-valued function $\mu(\cdot)$ on ${\bf X}$ is said to be a uniform Lyapunov function if the following conditions are satisfied:
\begin{itemize}
\item[(a)] $1+\sum_{y\in{\bf X}\setminus\{0\}} p(y|x,a)\mu(y)\le \mu(x),~\forall~x\in{\bf X},a\in{\bf A}.$
\item[(b)] For each $x\in{\bf X}$, the mapping $a\in {\bf A}\rightarrow\sum_{y\in{\bf X}\setminus\{0\}} p(y|x,a)\mu(y)$ is continuous.
\item[(c)] For each $x\in{\bf X}$ and each deterministic stationary strategy $\varphi$, $\lim_{t\to\infty} \EE^\varphi_x[\mu(X_t)$ $\times\I\{\tau_0>t\}]=0$, where we recall that $\tau_0:=\min\{t\ge 1:~X_t=0\}$.
\end{itemize}
\end{definition}
The above definition of a uniform Lyapunov function was taken from \cite[Definition 7.4]{b13}, see also \cite[Definition 4.2]{Cavazos:1992}. In both of these two references, this definition is ascribed to \cite{Hordijk:1975}.  Many characterizations and consequences of a semi-continuous MDP model with a uniform Lyapunov function can be found in \cite{b13,Cavazos:1992}, among which is the following one, whose proof can be found on p.107 of \cite{b13}.
\begin{proposition}\label{t2_04}
Consider a semi-continuous MDP model with a countable or finite state space $\bf X$ (endowed with the discrete topology) with the isolated absorbing state $0$. If there exists a uniform Lyapunov function $\mu(\cdot)$, then the MDP model is uniformly absorbing (at $0$) for each initial state $x_0.$
\end{proposition}

For a semi-continuous MDP model with a countable or finite state space $\bf X$ (endowed with the discrete topology) with the isolated absorbing state $0$,  if there exists a uniform Lyapunov function $\mu(\cdot)$, it can happen that for some initial distribution $\PP_0$, the MDP model is not absorbing (at $0$) for $\PP_0$. Sometimes, this is overlooked. We present an example to illustrate this.
\begin{example}\label{PiunovskiyZhang2023cExample02}
Consider an MDP model $\{{\bf X}, {\bf A},p\}$ with ${\bf X}=\{0,1,\dots\}$, endowed with the discrete topology, and ${\bf A}$ being a singleton, so that we will omit the argument $a\in{\bf A}$ everywhere in this example, and $p(0|x)=\frac{1}{2^x}=1-p(x|x)$ for all $x\in \{1,2,\dots\}$ and $p(0|0)=1$. There is only one strategy say $\pi$ in this model. Then for each $x_0\in\{1,2,\dots\}$,
\begin{eqnarray*}
\EE_{x_0}^\pi[T_0]=2^{x_0}.
\end{eqnarray*}
Hence, this MDP model is absorbing, in fact uniformly absorbing (at $0$) for every initial state $x_0$. It has a uniform Lyapunov function given by $\mu(0)=1$ and $\mu(x)=2^x$ for $x\ge 1.$ Indeed, (a,b) in Definition \ref{d2_01} can be checked to be satisfied by $\mu(\cdot)$, and for (c), note for each $x\ge 1$ that
\begin{eqnarray*}
\EE_x^\pi[\mu(X_n)\I\{\tau_0>n\}]=\EE_{x}^\pi[2^x \I\{X_n=x\}]=2^x \PP_x^\pi(X_n=x)=2^x (1-\frac{1}{2^x})^{n}\rightarrow 0
\end{eqnarray*}
as $n\rightarrow \infty$, and the convergence also takes place for $x=0$.
Now if we take the initial distribution as the geometric distribution on $\{1,2,\dots\}$ with parameter $\frac{2}{5}$, then $\EE_{\PP_0}[T_0]=\sum_{i=1}^\infty \frac{2}{5} (\frac{3}{5})^{i-1} 2^i=\infty.$
\end{example}

\subsubsection{Projection mapping and the topologies on ${\cal D}$ and ${\cal P}$}
Consider an absorbing (at $0$) MDP model for the initial distribution $\PP_0$. Let the projection mapping $O$ from ${\cal P}$ to ${\cal D}$ be defined by
\begin{eqnarray}\label{PiunovskiyZhang2023cEqn01}
O(\PP^\pi_{\PP_0}):=\eta_{\PP_0}^\pi.
\end{eqnarray}
The main results of this paper concern the continuity of the projection mapping $O$ from ${\cal P}$ to ${\cal D}$. Therefore, we endow ${\cal D}$ and ${\cal P}$ with suitable topologies, described as follows.

Consider an absorbing (at $0$) MDP model for the initial distribution $\PP_0$. Then ${\cal D}$ is a subset of the space ${\cal M}(({\bf X}\setminus\{0\}))\times {\bf A})$ of finite measures on $((({\bf X}\setminus\{0\})\times {\bf A},{\cal B}(({\bf X}\setminus\{0\})\times {\bf A}))$. We endow ${\cal M}(({\bf X}\setminus\{0\})\times {\bf A})$ with the weak topology generated by the set of bounded continuous functions defined on ${\bf X}\setminus\{0\})\times {\bf A}$, i.e., the coarsest topology on ${\cal M}((({\bf X}\setminus\{0\})\times {\bf A})$ with respect to which, for each bounded continuous function $g(\cdot)$ on $({\bf X}\setminus\{0\})\times {\bf A}$, $\int_{({\bf X}\times \{0\})\times {\bf A}}g(x,a)\eta(dx\times da)$ is continuous in $\eta\in {\cal M}(({\bf X}\setminus\{0\})\times {\bf A}).$  Since $\textbf{X}\setminus\{0\}$ and $\textbf{A}$ are Borel spaces, according to Theorem 8.9.4(i) of \cite{Bogachev:2007} and \cite[p.49]{Bogachev:2018}, when endowed the weak topology,  ${\cal M}(({\bf X}\setminus\{0\})\times {\bf A})$ is metrrizable, and its topological subspace ${\cal D}$ is  metrizable, too. Indeed, a popular topology on ${\cal D}$ is this weak topology, when ${\cal D}\subseteq {\cal M}(({\bf X}\setminus\{0\})\times {\bf A}).$

On the space ${\cal P}$ (for the fixed initial distribution $\PP_0$), one can consider the so called ws$^\infty$-topology of Sch\"al, introduced in \cite{b18}, see p.359 therein. This is the coarsest topology on ${\cal P}$, with respect to which, for each integer $0\le T<\infty$ and each bounded measurable function $f(h_T)=f(x_0,a_1,x_1,\ldots,a_T)$  continuous in $(a_1,a_2,\ldots,a_T)$ under arbitrarily fixed $x_0,x_1,\ldots,x_{T-1}\in {\bf X}$, the mapping $\PP\in {\cal P}\to\int_{{\bf H}_T} f(h_T)\PP_T(d h_T)$ is continuous. Here ${\bf H}_0={\bf X}$, and ${\bf H}_T=({\bf X}\times {\bf A})^{T}\times{\bf X}$, and ${\PP_T}$ denotes the marginal of $\PP\in {\cal P}$ on ${\bf H}_{T}$, $T\ge 0$ being an integer.

The ws$^\infty$-topology works particularly well with semi-continuous (S) MDP models. One of the main results in \cite{b18} is the following compactness result, whereas for semi-continuous (W) MDP models, the weak topology on ${\cal P}$ is more convenient. What actually happens is that for semi-continuous (S) MDP models, on ${\cal P}$, the ws$^{\infty}$-topology is the same as the weak topology, see Proposition \ref{PiunovskiyZhang2023cPropositionNowak} below.
\begin{proposition} \label{t009}
Consider the MDP model $\{{\bf X},{\bf A},p\}$ with some initial distribution $\PP_0$. If the MDP model is semi-continuous (S), then ${\cal P}$ endowed with ws$^{\infty}$ topology is compact, whereas if the MDP model is semi-continuous (W), then ${\cal P}$ endowed with the weak topology is compact.
\end{proposition}
\par\noindent\textit{Proof}. See \cite[Theorem 6.6]{b18} for the first assertion, and \cite[Theorem 5.6]{b18} for the second assertion. $\hfill\Box$

Nowak further studied the ws$^{\infty}$-topology on  ${\cal P}$ in \cite{Nowak:1988}, and proved the following useful result.
\begin{proposition}\label{PiunovskiyZhang2023cPropositionNowak}
Suppose that the MDP model $\{{\bf X},{\bf A},p\}$ is semi-continuous (S). Let some initial distribution $\PP_0$ be given. Then on ${\cal P}$, the ws$^{\infty}$-topology coincides with the weak topology. Consequently, the space ${\cal P}$ endowed with the ws$^\infty$-topology is metrizable, and compact in the weak topology, in view of Proposition \ref{t009}.
\end{proposition}
\par\noindent\textit{Proof.} This follows from \cite[Theorem 1]{Nowak:1988}.  A simpler proof of \cite[Theorem 1]{Nowak:1988} was given in \cite{Kurushima:2018}. $\hfill\Box$

In view of Propositions \ref{t009} and \ref{PiunovskiyZhang2023cPropositionNowak}, we will always consider the weak topology on ${\cal P}$ when we discuss semi-continuous MDP models, for which, ${\cal P}$ is always compact.

\subsection{Facts about optimal control problems}

Now we add an additional element to the MDP model, which is the cost function $c(\cdot)$. This is a measurable function from  $({\bf X}\times{\bf A},{\cal B}({\bf X}\times {\bf A}))$ to $[-\infty,\infty]$ in general.

Define for each strategy $\pi$ and an initial distribution $\PP_0$
 \begin{eqnarray*}
 v^\pi(\PP_0):=\EE^\pi_{P_0}\left[\sum_{t=1}^\infty c(X_{t-1},A_t)\right]:=\EE^\pi_{P_0}\left[\sum_{t=1}^\infty c^+(X_{t-1},A_t)\right]-\EE^\pi_{P_0}\left[\sum_{t=1}^\infty c^-(X_{t-1},A_t)\right],
 \end{eqnarray*}
 where $c^+(x,a):=\max\{c(x,a),0\}$, $c^-(x,a)=\max\{-c(x,a),0\}$, and we accept that $\infty-\infty:=\infty.$
When $\PP_0(dx)=\delta_{x_0}(dx)$ for some $x_0\in{\bf X}$, then we write $v^\pi(x_0)$ instead of $v^\pi(\PP_0).$

For future reference, let
\begin{eqnarray}\label{PiunovskiyZhang2023cEqn17}
v^\ast(x):=\inf_{\pi\in \Delta^{\rm All}}v^\pi(x)
\end{eqnarray}

A strategy $\pi^*$ is called optimal for $\PP_0$ if it  solves the following optimal control problem:
\begin{eqnarray}\label{PiunovskiyZhang2023cEqn11}
\mbox{Minimize over $\pi\in\Delta^{\rm All}$:} v^\pi(\PP_0),
\end{eqnarray} i.e.,
\begin{eqnarray*}
v^{\pi^*}(\PP_0)=\min_{\pi\in\Delta^{\rm All}}v^\pi(\PP_0).
\end{eqnarray*}

In the next statement, sufficient conditions for the existence of an optimal strategy for $\PP_0$ are given. See also \cite{Dufour:2020} and the references therein.
\begin{proposition} \label{t0091}
Suppose that the MDP model is semi-continuous (S), and $c(\cdot)$ is bounded below and $(-\infty,\infty]$-valued on ${\bf X}\times {\bf A}$ such that $c(x,a)$ is lower semicontinous in  $a\in{\bf A}$ for each $x\in{\bf X}$. If for the given initial distribution $\PP_0$,
 \begin{eqnarray}\label{PiunovskiyZhang2023cEqn12}
 \EE^\pi_{\PP_0}\left[\sum_{t=1}^\infty c^-(X_{t-1},A_t)\right]<\infty~\forall~\pi\in\Delta^{\rm All},
 \end{eqnarray}
 and
\begin{equation}\label{e1}
\inf_{N\ge n}\inf_{\pi\in\Delta^{\rm All}}\sum_{t=n+1}^N \EE^\pi_{\PP_0}\left[c(X_{t-1},A_t)\right]\to 0~\mbox{ as } n\to\infty,
\end{equation}
then there exists an optimal strategy for $\PP_0$. The same assertion holds if the MDP model is semi-continuous (W), and $c(\cdot)$ is a bounded below $(-\infty,\infty]$-valued lower semicontinous function on ${\bf X}\times {\bf A}$.
\end{proposition}
\par\noindent\textit{Proof.} The proof of the first assertion follows from \cite[Theorems 2.2 and 6.6]{b18}. The last assertion holds by \cite[Theorems 2.2 and 5.6]{b18} $\hfill\Box$

In \cite{b18}, condition (\ref{e1}) is called (C), and condition (\ref{PiunovskiyZhang2023cEqn12}) is referred to as ``General Assumption'' (GA).

For an absorbing (at $0$) MDP model for $\PP_0$,  we will take $c(\cdot)$ such that $c(0,a)\equiv 0$ at the absorbing isolated point $0$. Then $0$ will be referred to as the costless cemetery. If Condition (GA) is satisfied, then the optimal control problem (\ref{PiunovskiyZhang2023cEqn11}) can be rewritten in terms of occupation measures in the following way:
  \begin{eqnarray*}\label{e2_100}
\mbox{Minimize over $\eta\in {\cal D}$:~} \int_{{\bf X\setminus\{0\}}\times{\bf A}} c(y,a)\eta(dy\times da),
  \end{eqnarray*}
This is one of the main reasons for studying ${\cal D}$ in the literature of MDPs.

\section{Main results}\label{PiunovskiyZhang2023cSectionMain}

In the wonderful and insightful paper \cite{b114}, a rich theory for absorbing MDP models was developed. In that paper, the following assertion was claimed, see the proof of \cite[Lemma 4.7]{b114}:
\begin{itemize}
\item If the MDP model $\{{\bf X},{\bf A},p\}$ is semi-continuous, and is absorbing (at $0$) for the given initial distribution $\PP_0,$ ${\cal P}$ is endowed with the weak topology (see Proposition \ref{PiunovskiyZhang2023cPropositionNowak} and the paragraph below it), and $\cal D$ is endowed with the weak topology generated by bounded continuous functions on $({\bf X}\setminus\{0\})\times {\bf A}$, then the projection mapping $O$ defined in (\ref{PiunovskiyZhang2023cEqn01}) is continuous.
\end{itemize}
If the above claim was true, then by Propositions \ref{t009} and \ref{PiunovskiyZhang2023cPropositionNowak} as well as the fact ${\cal D}= O({\cal P})$, it would then follow that ${\cal D}$ is compact, provided that the MDP model $\{{\bf X},{\bf A},p\}$ is semi-continuous, and is absorbing (at $0$) for $\PP_0.$ This assertion was formulated as \cite[Lemma 4.7]{b114}, and was used several times therein. 

The first main result of this paper is that we show by means of an example that the above claims regarding the continuity of $O$ and the compactness of ${\cal D}$ are false.
\begin{theorem}\label{PiunovskiyZhang2023cThm01}
There is an MDP model $\{{\bf X},{\bf A},p\}$ such that the following assertions hold for a fixed initial distribution $\PP_0$:
\begin{itemize}
\item[(a)] The MDP model is semi-continuous. Consequently, the space of all strategic measures ${\cal P}$ for the initial distribution $\PP_0$ is compact in the weak topology.
\item[(b)] The MDP model is absorbing (at $0$) for $\PP_0.$
\item[(c)] The MDP model is not uniformly absorbing (at $0$) for $\PP_0$.
\item[(d)] ${\cal D}$ is not compact with respect to the weak topology.
\item[(e)] The mapping $O$ defined by (\ref{PiunovskiyZhang2023cEqn01}) is not continuous between ${\cal P}$ and ${\cal D}$, both of which are endowed with their weak topologies.
\item[(f)] The MDP model does not possess a uniform Lyapunov function.
\end{itemize}
\end{theorem}
The proof of this theorem is given in Subsection \ref{PiunovskiyZhang2023cSubsection01}.

Since the continuity of the projection mapping $O$ defined by (\ref{PiunovskiyZhang2023cEqn01}) from ${\cal P}$ (endowed with the weak topology) to  ${\cal D}$ (endowed with the weak topology) played an important role in the reasoning in \cite{b114}, and Theorem \ref{PiunovskiyZhang2023cThm01} shows that this is not guaranteed if the MDP model is semi-continuous and absorbing (at $0$) for $\PP_0$, we next investigate when it holds. Our second main result asserts that for a semi-continuous MDP model, provided that it is absorbing for $\PP_0$, the aforementioned continuity of the projection mapping $O$ holds if and only if the MDP model is uniformly absorbing for $\PP_0$.

\begin{theorem}\label{t92}
Consider a semi-continuous MDP model with the fixed initial distribution $\PP_0$. Suppose that the MDP model is absorbing (at $0$) for $\PP_0$. Then the projection mapping $O$ defined by (\ref{PiunovskiyZhang2023cEqn01}) from ${\cal P}$ (endowed with the weak topology) to ${\cal D}$ (endowed with the weak topology) is continuous if and only if the MDP model is uniformly absorbing (at $0$) for $\PP_0$.
\end{theorem}
The proof is given in Subsection \ref{PiunovskiyZhang2023cSubsection02}.

\section{Example}\label{PiunovskiyZhang2023cSectionExample}

Theorem \ref{PiunovskiyZhang2023cThm01}  will be proved by means of the following MDP model.

\begin{figure}[ht]
\begin{center}
\includegraphics[height=6.5cm]{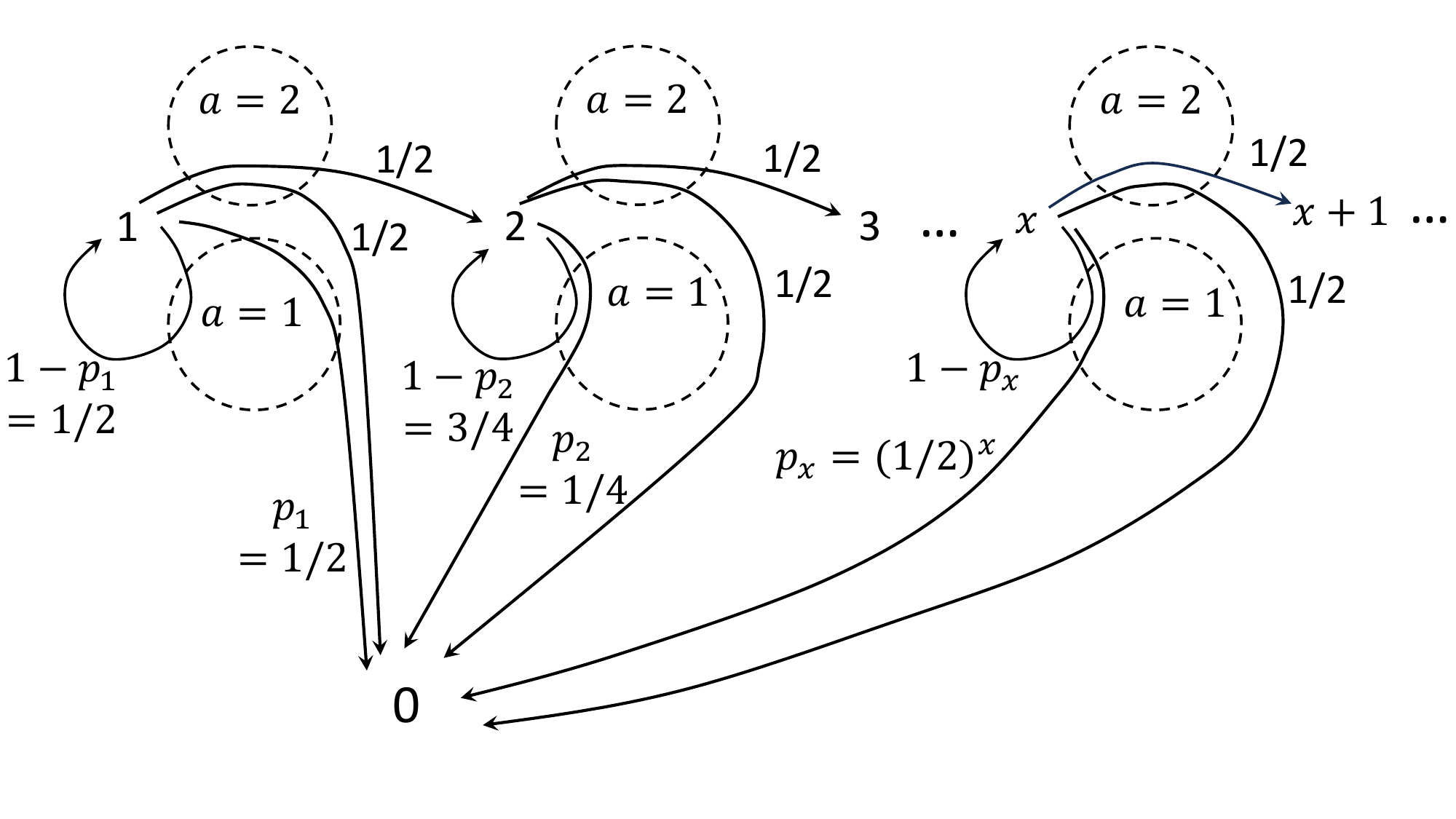}
\end{center}
\caption{\label{f2_17} Graphical representation of the MDP.}
\end{figure}

\begin{example}\label{PiunovskiyZhang2023cExample01}
The elements of the MDP model $\{{\bf X},{\bf A},p\}$ are as follows.
\begin{itemize}
\item The state and action spaces are ${\bf X}:=\{0,1,2,\ldots\}$ and ${\bf A}:=\{1,2\}$, which are equipped with the discrete topology.
\item $0$ is the isolated state, which is absorbing: $p(0|0,a)\equiv 0$.
\item For $x=1,2,\ldots$, $p(0|x,1)=p_x=\left(\frac{1}{2}\right)^x$ and $p(x|x,1)=1-p_x$.
\item For $x=1,2,\ldots$, $p(0|x,2)=\frac{1}{2}$ and $p(x+1|x,2)=\frac{1}{2}$.
\item Other transition probabilities equal zero.
\end{itemize}
\end{example}
The transition diagram of the MDP model in Example \ref{PiunovskiyZhang2023cExample01} is given in Figure \ref{f2_17}.

The MDP model in Example \ref{PiunovskiyZhang2023cExample01} is explored in the proof of Theorem \ref{PiunovskiyZhang2023cThm01}. Incidentally, it can also be used to demonstrate that condition (\ref{e1}) is important for the solvability, see the next proposition.

\begin{proposition}\label{PiunovskiyZhang2023cThm02}
Consider the MDP model in Example \ref{PiunovskiyZhang2023cExample01} with the initial distribution \linebreak $\PP_0(dx)=\delta_1(dx)$, and the optimal control problem (\ref{PiunovskiyZhang2023cEqn11}) for the cost function $c(\cdot)$ defined by $c(x,a)\equiv -1$ for all $x\ge 1$, and $c(0,a)\equiv 0.$ Then the following assertions hold.
\begin{itemize}
\item[(a)] There is no optimal strategy for the given $\PP_0$.
\item[(b)] Condition (\ref{e1}) is not satisfied, whereas all the other conditions in Proposition \ref{t0091} are satisfied.
\end{itemize}
\end{proposition}
The proof of this Proposition is given in Subsection \ref{PiunovskiyZhang2023cSectionTheorme3}.

\section{Proof of the statements}\label{PiunovskiyZhang2023cSectionProof}

In this section we provide the detailed proofs of Theorem \ref{PiunovskiyZhang2023cThm01}, Theorem \ref{t92} and Proposition \ref{PiunovskiyZhang2023cThm02}.

\subsection{Proof of Theorem \ref{PiunovskiyZhang2023cThm01}}\label{PiunovskiyZhang2023cSubsection01}
\par\noindent\textit{Proof of Theorem \ref{PiunovskiyZhang2023cThm01}.} Throughout this proof we consider the MDP model in Example \ref{PiunovskiyZhang2023cExample01}. Unless stated otherwise, we fix the initial distribution $\PP_0$ concentrated on the singleton $\{1\}$, i.e.,
\begin{eqnarray*}
\PP_0(1)=1.
\end{eqnarray*}

(a) This MDP is vacuously semi-continuous, and the space of all strategic measures ${\cal P}=\{\PP^\pi_{\PP_0},~\pi\in\Delta^{\rm All}\}$ is compact in the weak topology by Propositions \ref{t009} and \ref{PiunovskiyZhang2023cPropositionNowak}.

(b) Consider the function $w^\ast(\cdot)$ on ${\bf X}=\{0,1,\dots,\}$ defined by
\begin{eqnarray}\label{PiunovskiyZhang2023cEqn06}
w^\ast(x)&:=&\sup_{\pi\in\Delta^{\rm All}}\EE_{x}^\pi[T_0]=\sup_{\pi\in\Delta^{\rm All}}\EE_{x}^\pi\left[\sum_{t=0}^\infty \I\{X_t\in\{1,2,\dots\}\}\right]\\
&=&\sup_{\pi\in\Delta^{\rm All}}\sum_{y=1}^\infty\eta^\pi_x(y),
~x\in\{0,1,\dots\}.\nonumber
\end{eqnarray}
We shall show that $w^\ast(x)<\infty$ for each $x\in \{0,1,\dots\}$. This would in particular, imply that $\EE_{1}^\pi[T_0]<\infty$ for all $\pi\in\Delta^{\rm All}$, i.e., the MDP model is absorbing (at $0$) for the fixed $P_0$ concentrated on $\{1\}$.

According to \cite[Propositions 9.8 and 9.10]{b7},
$w^\ast(\cdot)$ is the minimal nonnegative solution to the Bellman equation
\begin{eqnarray}\label{PiunovskiyZhang2023cEqn02}
w^*(0) &=& 0;\nonumber\\
w^*(x) &=& 1 +\max\left\{(1-p_x)w^*(x),~\frac{1}{2} w^*(x+1)\right\},~x=1,2,\ldots.
\end{eqnarray}
We claim that $w^\ast(0)=0$ and
\begin{eqnarray}\label{PiunovskiyZhang2023cEqn07}
w^\ast(x)\le 2+2^x~\forall~x\in\{1,2,\dots\}.
\end{eqnarray} To see this, it suffices to note that
the function $v(\cdot)$ on $\{0,1,\dots\}$ given by $v(x):=2+2^x$ and $v(0)=0$ solves the Bellman equation (\ref{PiunovskiyZhang2023cEqn02}), and this is true because, $v(0)=0$ trivially satisfies the first equality in (\ref{PiunovskiyZhang2023cEqn02}), whereas for $x\in\{1,2,\dots\}$,
\begin{eqnarray*}
&& 1 +\max\left\{\left[1-\left(\frac{1}{2}\right)^x\right][2+2^x];~\frac{1}{2}(2+2^{x+1})\right\}\\
&=&  1+\max\left\{1+2^x-\left(\frac{1}{2}\right)^{x-1};~1+2^x\right\}=2+2^x.
\end{eqnarray*}
Therefore, $w^\ast(0)= 0$ and $w^\ast(x)\le 2+2^x$ for $x\in\{1,2,\dots\}$. Thus, (b) is proved.

Incidentally, from the calculations similar to those in the proof of (c), we will see actually that
\begin{eqnarray*}
w^\ast(x)= 2+2^x ~\forall~x\in\{1,2,\dots\}.
\end{eqnarray*}
This equality will be used in the proof of a subsequent statement, and we formulate it as Lemma \ref{PiunovskiyZhang2023cLemma01} below, and will prove it there. Nevertheless, for the purpose here, the validity of the inequality $\le$ is sufficient.

(c) We have seen in part (b) that this MDP model is absorbing (at $0$) for the given initial distribution $\PP_0$ concentrated on $\{1\}$. We now verify that it is not uniformly absorbing (at $0$) for $\PP_0$.

To show this, consider the deterministic stationary strategies
\begin{eqnarray}\label{PiunovskiyZhang2023cEqn08}
\varphi^n(x):= 2\cdot\I\{x\le n\}+\I\{x>n\}~\forall~ n=0,1,2,\ldots.
\end{eqnarray}
With slight abuse of notation, let us denote by
\begin{eqnarray*}
\eta_{\PP_0}^{\pi}(x):=\eta_{\PP_0}^\pi(\{x\}\times{\bf A})
\end{eqnarray*}
the marginal on ${\bf X}\setminus \{0\}=\{1,2,\dots\}$ of the occupation measure of the strategy $\pi$ for the initial distribution $\PP_0$.

For each $n\in\{0,1,\dots\}$, $\eta_1^{\varphi^n}$ is given by
\begin{equation}\label{e55}
\eta_1^{\varphi^n}(x)=\left\{\begin{array}{ll}
\left(\frac{1}{2}\right)^{x-1}, & \mbox{ if } x<n+1;\\
\left(\frac{1}{2}\right)^{n}\frac{1}{p_{n+1}}=2, & \mbox{ if } x=n+1;\\
0, & \mbox{ if } x>n+1.\end{array}\right.
\end{equation}
Indeed, under the deterministic stationary strategy $\varphi^n$ any state $x<n+1$ is reached with probability $\left(\frac{1}{2}\right)^{x-1}$, and given that it is reached, the controlled process spends exactly one time unit on it; any state $x>n+1$ is never reached. This justifies the first and the third equalities in (\ref{e55}). For the state $x=n+1$, note that $\left(\frac{1}{2}\right)^n$ is the probability that $X_n=n+1.$ (The complementary probability $1-\left(\frac{1}{2}\right)^n$ is for $X_n=0$; the state $n+1$ cannot appear at the steps $t<n$.) After that $X_{n+t}=n+1$ is realized with probability $(1-p_{n+1})^t$, leading to
\begin{eqnarray}\label{e33}
\eta_1^{\varphi^n}(n+1)=\left(\frac{1}{2}\right)^n\left[1+(1-p_{n+1})+(1-p_{n+1})^2-\ldots\right]=\left(\frac{1}{2}\right)^n\frac{1}{p_{n+1}}=2.
\end{eqnarray}
For each $n\ge 0$,
\begin{eqnarray}\label{PiunovskiyZhang2023cEqn16}
\sup_{\pi\in\Delta^{\rm All}} \EE^\pi_{1}\left[\sum_{t=n}^\infty \I\{X_t\in({\bf X}\setminus\{0\})\}\right]\ge \EE^{\varphi^n}_{1}\left[\sum_{t=n}^\infty \I\{X_t=n+1\}\right]=\eta^{\varphi^n}_1(n+1)=2,
\end{eqnarray}
so that $\sup_{\pi\in\Delta^{\rm All}} \EE^\pi_{1}\left[\sum_{t=n}^\infty \I\{X_t\in({\bf X}\setminus\{0\})\}\right]$ does not converge to $0$ as $n\rightarrow \infty$. In view of (\ref{PiunovskiyZhang2023cEqn05}), it follows that this MDP model is not uniformly absorbing at $0$ for the given $\PP_0$.

(d) Now we fix the weak topology on ${\cal D}$, as described in Section \ref{PiunovskiyZhang2023cSectionFacts}. Accordingly, in this proof, the notions of compactness and convergence of sequences in ${\cal D}$ are understood with respect to it, and this will not be signified repeatedly. The target here is to show that the space $\cal D$ is not compact.  This is equivalent to showing that ${\cal D}$ is not sequentially compact, because the weak topology on ${\cal D}$ is metrizable, as mentioned in Section \ref{PiunovskiyZhang2023cSectionFacts}. Consequently, it suffices to show that the set $\{\eta_1^{\varphi^n},~n=0,1,2,\ldots\}$ of occupation measures on ${\cal B}(({\bf X}\setminus\{0\})\times{\bf A})$ has no accumulation points in ${\cal D}$, where the deterministic stationary strategies $\varphi^n$ are defined by (\ref{PiunovskiyZhang2023cEqn08}) in the proof of (c).

Suppose for contradiction that $\eta_1^{\varphi^{n_i}}\to \eta_1^\pi\in {\cal D}$ for some subsequence $n_i\to\infty$ as $i\to\infty$, and for some strategy $\pi$.
For a fixed $j=1,2,\ldots$ take $d^j(x,a):=\I\{x=j\}$. Then
\begin{eqnarray}\label{PiunovskiyZhang2023cEqn10}
\lim_{i\to\infty} \sum_{(x,a)\in({\bf X}\setminus\{0\})\times{\bf A}} d^j(x,a)\eta_1^{\varphi^{n_i}}(x,a)=\left(\frac{1}{2}\right)^{j-1}
\end{eqnarray}
according to (\ref{e55}), because $n_i\rightarrow \infty$ as $i\rightarrow \infty.$

Since $d^j(\cdot)$ is bounded continuous on $({\bf X}\setminus\{0\})\times {\bf A}$, from $\eta_1^{\varphi^{n_i}}\to \eta_1^\pi\in {\cal D}$ we see that
\begin{eqnarray}\label{PiunovskiyZhang2023cEqn09}
\lim_{i\to\infty} \sum_{(x,a)\in({\bf X}\setminus\{0\})\times{\bf A}} d^j(x,a)\eta_1^{\varphi^{n_i}}(x,a)= \sum_{(x,a)\in({\bf X}\setminus\{0\})\times{\bf A}} d^j(x,a)\eta_1^{\pi}(x,a),~\forall~j\ge 1.
\end{eqnarray}

Applying (\ref{PiunovskiyZhang2023cEqn10}) and (\ref{PiunovskiyZhang2023cEqn09}) to $j=1$, we see that
\begin{eqnarray*}
\sum_{(x,a)\in({\bf X}\setminus\{0\})\times{\bf A}} d^1(x,a)\eta_1^{\pi}(x,a)=\left(\frac{1}{2}\right)^0=1.
\end{eqnarray*}
On the other hand, if $\pi_1(1|1)=\varepsilon\ge 0$ and $\pi_1(2|1)=1-\varepsilon$, then for $j=1$, we have
\begin{eqnarray*}
\sum_{(x,a)\in({\bf X}\setminus\{0\})\times{\bf A}} d^1(x,a)\eta_1^{\pi}(x,a)\ge
\EE^\pi_1\left[\I\{X_0=1\}+\I\{X_1=1\}\right]=
 1+\varepsilon(1-p_1),
\end{eqnarray*}
so that, it is necessary that $\varepsilon=0$.

The similar calculation for $j=2$, $d^2(x,a)=\I\{x=2\}$, and $\pi_2(a|x_0,a_1,x_1)$ with $x_0=1,a_1=2$ and $x_1=2$ leads to equality $\pi_2(2|1,2,2)=1$, and so on. Thus, $\pi$ has the same occupation measure as the deterministic stationary strategy given by $\varphi(x)\equiv 2$, i.e., $\eta_1^\pi=\eta_1^{\varphi}.$

 However, for the bounded continuous function $d(\cdot)$ on $({\bf X}\setminus \{0\})\times {\bf A}$ given by $d(x,a)\equiv 1$,
\begin{eqnarray}\label{PiunovskiyZhang2023cEqn18}
\sum_{(x,a)\in({\bf X}\setminus\{0\})\times{\bf A}} d(x,a)\eta_1^{\varphi}(x,a)=1+\frac{1}{2}+\frac{1}{4}+\ldots=2
\end{eqnarray}
whereas
\begin{eqnarray}\label{PiunovskiyZhang2023cEqn19}
\sum_{(x,a)\in({\bf X}\setminus\{0\})\times{\bf A}} d(x,a)\eta_1^{\varphi^{n_i}}(x,a)=\sum_{x=1}^{n_i} \left(\frac{1}{2}\right)^{x-1}+2=4-\frac{1}{2^{n_i-1}}\to 4 ~\mbox{ as } i\to\infty.
\end{eqnarray}
The previous two equalities yield the desired contradiction against the assumption $\eta_1^{\varphi^{n_i}}\to \eta_1^\pi$, because $d(\cdot)$ is bounded and continuous on $({\bf X}\setminus \{0\})\times {\bf A}$. Thus, the space ${\cal D}$  is not compact with respect to the weak topology, as required.

(e) As mentioned in Section \ref{PiunovskiyZhang2023cSectionMain}, if the projection mapping $O$ defined by (\ref{PiunovskiyZhang2023cEqn01}) was continuous from ${\cal P}$ endowed with the weak topology to ${\cal D}$ endowed with the weak topology, then ${\cal D}$ would have been compact, because ${\cal P}$ is compact by Proposition \ref{t009}. Since ${\cal D}$ is not compact as shown in (d), $O$ is not continuous.

Alternatively, we may deduce it more explicitly, as follows. Consider the strategies $\varphi^{n}$ and $\varphi$ as in the proof of (d). Let us show that $\PP^{\varphi^n}_1\to \PP^\varphi_1$ as $n\to\infty$ in the weak topology. By Proposition \ref{PiunovskiyZhang2023cPropositionNowak}, it is sufficient to show that $\PP^{\varphi^n}_1\to \PP^\varphi_1$ as $n\to\infty$ in the ws$^{\infty}$-topology.
For a fixed $0\le T<\infty$, if $n\ge T$, then only the sequences $h_T=(x_0,a_1,x_1,\ldots,x_{T-1})$ of the form
\begin{eqnarray*}
&& (1,2,0,2,0,2,0,2,0,\ldots,0)\\
&& (1,2,2,2,0,2,0,2,0,\ldots,0)\\
&& (1,2,2,2,3,2,0,2,0,\ldots,0)\\
&& (1,2,2,2,3,2,4,2,0,\ldots,0)\\
&& (1,2,2,2,3,2,4,2,5,\ldots,0)\\
&& \ldots~~~\ldots\\
&& (1,2,2,2,3,2,4,2,5,\ldots,0)\\
&& (1,2,2,2,3,2,4,2,5,\ldots,T)
\end{eqnarray*}
can be realized with positive $\PP^{\varphi^n}_1$-probabilities $\frac{1}{2}$, $\frac{1}{4}$, $\frac{1}{8}$, $\frac{1}{16}$, $\frac{1}{32}$, $\ldots$, $\left(\frac{1}{2}\right)^{2T-3}$, $\left(\frac{1}{2}\right)^{2T-3}$, correspondingly. (Note, the last two probabilities coincide.) The same holds true for $\PP^\varphi_1$. Consequently, for each $T\ge 0,$
\begin{eqnarray*}
\lim_{n\to\infty} \int_{\bf H} f(h_T)d\PP^{\varphi^n}_1(h_T)= \int_{\bf H} f(h_T)d\PP^{\varphi}_1(h_T)
\end{eqnarray*}
for each bounded function $f(h_T)$, which is certainly continuous in the discrete topology on $\bf X$ and $\bf A$. Therefore, $\PP^{\varphi^n}_1\to \PP^\varphi_1$ as $n\to\infty$ in the ws$^\infty$-topology, and thus in the weak topology.

However, in the proof of (d), we have seen that $O(\PP_1^{\varphi^{n}})=\eta_1^{\varphi^{n}}$ does not converge to $O(\PP_1^\varphi)=\eta_1^\varphi$ in the weak topology, see (\ref{PiunovskiyZhang2023cEqn18}) and (\ref{PiunovskiyZhang2023cEqn19}). This shows the claimed discontinuity of the mapping $O$.

(f) Since this MDP model is not uniformly absorbing (at $0$) for the initial state $x_0=1$, it follows from Proposition \ref{t2_04} that there cannot be uniform Lyapunov functions. In fact, conditions (a,b) in Definition \ref{d2_01} can be satisfied, but any function satisfying them violates condition (c) in Definition \ref{d2_01}. We demonstrate this fact explicitly as follows. First, the function defined by $\mu(0)=1$ and $\mu(x)=2+2^x$ for $x\ge 1$ satisfies condition (a) in Definition \ref{d2_01} because
$1+\sum_{y=1}^\infty \mu(y)p(y|0,a)\equiv 1=\mu(0)$, and
 for each $x\ge 1,$ it holds that
\begin{eqnarray*}
&&1+\sum_{y=1}^\infty \mu(y)p(y|x,1)=1+(1-p_x)\mu(x)\\
&=&1+(1-\frac{1}{2^x})(2+2^x)=1+2+2^x-\frac{1}{2^{x-1}}-1
=2+2^x-\frac{1}{2^{x-1}}\le 2+2^x=\mu(x)\end{eqnarray*}
and
\begin{eqnarray*}
&&1+\sum_{y=1}^\infty \mu(y)p(y|x,2)=1+\frac{1}{2}\mu(x+1)=1+1+2^x=2+2^x=\mu(x).
\end{eqnarray*}
Condition (b) in Definition \ref{d2_01} is trivially satisfied.

Next, let us show that, for each function $\mu(\cdot)$ satisfying condition (a) in Definition \ref{d2_01}, condition (c) therein is violated. Indeed, condition (a) in Definition \ref{d2_01} implies that
\begin{eqnarray*}
1+\max\left\{(1-p_x)\mu(x);~\frac{1}{2}\mu(x+1)\right\}\le \mu(x),~\forall~x>0,
\end{eqnarray*}
leading to, for each $x\ge 1$,
\begin{eqnarray*}
&&1+(1-p_x)\mu(x)\le \mu(x)\Longrightarrow \mu(x)\ge \frac{1}{p_x}=2^x;\\
&&\mu(x)\ge 1+\frac{1}{2}\mu(x+1)\ge 1+\frac{1}{2}\cdot 2^{x+1}=1+2^x.
\end{eqnarray*}
Now, for the deterministic stationary strategy $\varphi(x)\equiv 2$, we have for each $x\ge 1$ that
\begin{eqnarray*}
\EE^\varphi_x[\mu(X_t)\cdot\I\{T_0>t\}]\ge \left(\frac{1}{2}\right)^t[1+2^{x+t}]=\left(\frac{1}{2}\right)^t+2^x,
\end{eqnarray*}
where $\left(\frac{1}{2}\right)^t=\PP^\varphi_x(T_0>t)$ and, if $T_0>t$, then $X_t=x+t$ with probability $1$.
Since the above expression does not converge to $0$ as $t\rightarrow \infty$, we see that condition (c) in Definition \ref{d2_01} is not satisfied by $\mu(\cdot)$. $\hfill\Box$

\subsection{Proof of Theorem \ref{t92}}\label{PiunovskiyZhang2023cSubsection02}
\par\noindent\textit{Proof of Theorem \ref{t92}.}
Suppose the MDP model is semi-continuous. First, we prove the `if' part. This is done by mimicking the reasoning in the proof of \cite[Lemma 4.7]{b114}, avoiding the minor inaccuracy therein. Suppose $\PP^{\pi^i}_{\PP_0}\to \PP^{\hat\pi}_{\PP_0}$ as $i\to\infty$ in the weak topology, and fix an arbitrary $\varepsilon>0$ and an arbitrary bounded continuous function $g(\cdot)$ on ${\bf X}\times{\bf A}$ such that $g(0,a)\equiv 0$. Let
\begin{eqnarray*}
\bar g:=\sup_{(x,a)\in({\bf X}\setminus\{0\})\times{\bf A}} |g(x,a)|.
\end{eqnarray*}
Let $N\ge 0$ be such that, for all $\pi\in\Delta^{\rm All}$,
\begin{eqnarray*}
\bar g\sup_{\pi\in\Delta^{\rm All}} \EE^\pi_{\PP_0}\left[\sum_{t=N}^\infty\I\{X_t\in({\bf X}\setminus\{0\})\}\right]<\varepsilon/3.
\end{eqnarray*}
Such an $N\ge 0$ exists because the MDP model is assumed to be uniformly absorbing (at $0$) for $\PP_0$ here.
Then let $K$ be such that, for all $i\ge K$,
\begin{eqnarray*}
\left|\EE^{\pi^i}_{\PP_0}\left[\sum_{t=1}^{N-1} g(X_{t-1},A_t)\right]-\EE^{\hat\pi}_{\PP_0}\left[\sum_{t=1}^{N-1} g(X_{t-1},A_t)\right]\right|<\varepsilon/3.
\end{eqnarray*}
Such a $K$ exists because by assumption, $\PP^{\pi^i}_{\PP_0}\to \PP^{\hat\pi}_{\PP_0}$ as $i\to\infty$ in the weak topology.

Now, for all $i\ge K$,
\begin{eqnarray*}
&&\left|\int_{({\bf X}\setminus\{0\})\times{\bf A}} g(x,a) \eta_{\PP_0}^{\pi^i}(dx\times da)- \int_{({\bf X}\setminus\{0\})\times{\bf A}} g(x,a) \eta_{\PP_0}^{\hat\pi}(dx\times da)\right|\\
&=&\left|\EE^{\pi^i}_{\PP_0}\left[\sum_{t=1}^{\infty} g(X_{t-1},A_t)\right]-\EE^{\hat\pi}_{\PP_0}\left[\sum_{t=1}^{\infty} g(X_{t-1},A_t)\right]\right|\\
&\le &\left|\EE^{\pi^i}_{\PP_0}\left[\sum_{t=1}^{N-1} g(X_{t-1},A_t)\right]-\EE^{\hat\pi}_{\PP_0}\left[\sum_{t=1}^{N-1} g(X_{t-1},A_t)\right]\right|\\
&+& \EE^{\pi^i}_{\PP_0}\left[\sum_{t=N}^{\infty} |g(X_{t-1},A_t)|\right]+ \EE^{\hat\pi}_{\PP_0}\left[\sum_{t=N}^{\infty} |g(X_{t-1},A_t)|\right]<\varepsilon,
\end{eqnarray*}
implying that $\eta^{\pi^i}=O(\PP_{\PP_0}^{\pi^i})\to\eta^{\hat\pi}=O(\PP_{\PP_0}^{\hat{\pi}})$ as $i\to \infty$. Since the MDP model is semi-continuous, ${\cal P}$ is metrizable by Proposition \ref{PiunovskiyZhang2023cPropositionNowak}. It follows that the mapping $O$ defined by (\ref{PiunovskiyZhang2023cEqn01}) is continuous.

Second, we prove the `only if' part. Suppose the MDP model is semi-continuous and absorbing (at $0$) for $\PP_0$, but not uniformly absorbing (at $0$) for $\PP_0$. Then there is $\varepsilon>0$ such that, for each $k\ge 1$, there is $n>k$ such that $\sup_{\pi\in\Delta^{\rm All}} \EE^\pi_{\PP_0} \left[\sum_{t=n}^\infty \I\{X_t\ne 0\}\right]>\varepsilon$. Below, such $\varepsilon>0$ is fixed.

Now
\begin{itemize}
\item there is $n_1\ge 1$ such that there exists a strategy $\pi^1$ satisfying the inequality
\begin{eqnarray*}
\EE^{\pi_1}_{\PP_0} \left[\sum_{t=n_1}^\infty \I\{X_t\ne 0\}\right]>\varepsilon;
\end{eqnarray*}
\item there are $n_2>n_1$ and $\pi^2$ such that
\begin{eqnarray*}
\EE^{\pi_2}_{\PP_0} \left[\sum_{t=n_2}^\infty \I\{X_t\ne 0\}\right]>\varepsilon;
\end{eqnarray*}
\item and so on.
\end{itemize}
In this manner, we obtain sequences $\{n_i\}_{i=1}^\infty$ and $\{\pi^i\}_{i=1}^\infty$ such that $n_i\uparrow\infty$ and, for each $i=1,2,\ldots$,
\begin{eqnarray*}
\EE^{\pi^i}_{\PP_0} \left[\sum_{t=n_i}^\infty \I\{X_t\ne 0\}\right]>\varepsilon.
\end{eqnarray*}

Since $\cal D$ is compact and metrizable in the ws$^\infty$-topology by Propositions \ref{t009} and \ref{PiunovskiyZhang2023cPropositionNowak}, there exists a subsequence  $\{\pi^{i_j}\}_{j=1}^\infty$ (with the corresponding subsequence $\{n_{i_j}\}_{j=1}^\infty$, $\lim_{j\to\infty}n_{i_j}=\infty$) such that $\PP^{\pi^{i_j}}_{\PP_0}\to \PP^{\hat\pi}_{\PP_0}$ in the weak topology for some strategy $\hat\pi$. We will show that  the sequence $\{\eta_{\PP_0}^{\pi^{i_j}}\}_{j=1}^\infty$ does not converge to $\eta_{\PP_0}^{\hat\pi}$ in the weak topology.

For the constant $\varepsilon>0$ fixed above, we choose $N>0$ such that
\begin{eqnarray*}
\EE^{\hat\pi}_{\PP_0}\left[\sum_{t=N}^\infty \I\{X_{t-1}\ne 0\}\right]=\EE^{\hat\pi}_{\PP_0}\left[\sum_{t=N}^\infty \I\{t<T_0\}\right]<\varepsilon/3.
\end{eqnarray*}
This can be done because the MDP is absorbing (at $0$) for $\PP_0$: $\EE^{\hat\pi}_{\PP_0}[T_0]=\sum_{n=0}^\infty \PP^{\hat\pi}_{\PP_0}(T_0>n)<\infty$.

After that, choose $K>0$ such that, if $j>K$, then
\begin{eqnarray*}
\EE^{\pi^{i_j}}_{\PP_0}\left[\sum_{t=1}^{N-1} \I\{X_{t-1}\ne 0\}\right]-\EE^{\hat\pi}_{\PP_0}\left[\sum_{t=1}^{N-1} \I\{X_{t-1}\ne 0\}\right]>-\varepsilon/3.
\end{eqnarray*}
This can be done because $\PP^{\pi^{i_j}}_{\PP_0}\to \PP^{\hat\pi}_{\PP_0}$ in the weak topology. Recall that $0$ is an isolated state so that $\I\{x\ne 0\}$ is a bounded continuous function on ${\bf X}.$

Now, for $j>K$, $n_{i_j}>N$, i.e., for all big values of $j$, we have
\begin{eqnarray*}
\EE^{\pi^{i_j}}_{\PP_0} \left[\sum_{t=N}^\infty \I\{X_t\ne 0\}\right]\ge \EE^{\pi^{i_j}}_{\PP_0} \left[\sum_{t=n_{i_j}}^\infty \I\{X_t\ne 0\}\right]> \varepsilon
\end{eqnarray*}
and
the following relations for the bounded continuous function $g(x,a)\equiv 1$ on $({\bf X}\setminus\{0\})\times {\bf A}$:
\begin{eqnarray*}
&&\int_{({\bf X}\setminus\{0\})\times{\bf A}} g(x,a) \eta_{\PP_0}^{\pi^{i_j}}(dx\times da)- \int_{({\bf X}\setminus\{0\})\times{\bf A}} g(x,a) \eta^{\hat\pi}(dx\times da)\\
&=&\EE^{\pi^{i_j}}_{\PP_0}\left[\sum_{t=1}^{N-1} \I\{X_{t-1}\ne 0\}\right]-\EE^{\hat\pi}_{\PP_0}\left[\sum_{t=1}^{N-1} \I\{X_{t-1}\ne 0\}\right]\\
&&+ \EE^{\pi^{i_j}}_{\PP_0} \left[\sum_{t=N}^\infty \I\{X_t\ne 0\}\right]- \EE^{\hat\pi}_{\PP_0}\left[\sum_{t=N}^\infty \I\{X_{t-1}\ne 0\}\right]\\
&>&-\frac{\varepsilon}{3}+\varepsilon-\frac{\varepsilon}{3}=\varepsilon/3.
\end{eqnarray*}
Thus, the sequence $\{\eta_{\PP_0}^{\pi^{i_j}}\}_{j=1}^\infty$ does not converge to $\eta_{\PP_0}^{\hat\pi}$ in the weak topology, and the mapping $O$ defined by (\ref{PiunovskiyZhang2023cEqn01}) is not continuous. $\hfill\Box$

\subsection{Proof of Proposition \ref{PiunovskiyZhang2023cThm02}}\label{PiunovskiyZhang2023cSectionTheorme3}

In this subsection, we prove Proposition \ref{PiunovskiyZhang2023cThm02}. Firstly, we present a lemma.

\begin{lemma}\label{PiunovskiyZhang2023cLemma01}Consider the MDP model in Example \ref{PiunovskiyZhang2023cExample01}, and the function $w^\ast(\cdot)$ on ${\bf X}$ defined in (\ref{PiunovskiyZhang2023cEqn06}). Then
\begin{eqnarray*}
&&w^\ast(0)=0;\\
&&w^\ast(x)= 2+2^x ~\forall~x\in\{1,2,\dots\}.
\end{eqnarray*}
\end{lemma}

\par\noindent\textit{Proof.}The first equality holds vacuously. The rest verifies the second equality. For a fixed $x\ge 1$, consider the following deterministic stationary strategies
\begin{eqnarray*}
\varphi^{n}(y):= 2\cdot\I\{y\le x+n\}+\I\{y>x+n\},~~n=0,1,\ldots.
\end{eqnarray*}
Then similar considerations to those for (\ref{e55}) result in, with slight abuse of notation,
\begin{equation*}
\eta_x^{\varphi^{n}}(y)=\left\{\begin{array}{ll}
\left(\frac{1}{2}\right)^{y-x}, & \mbox{ if } x\le y < x+n+1;\\
\left(\frac{1}{2}\right)^{y-x}\frac{1}{p_{x+n+1}}=2^x, & \mbox{ if } y=x+n+1;\\
0, & \mbox{ if } y>x+n+1, \mbox{ or }y<x.\end{array}\right.
\end{equation*}
for the marginal on ${\bf X}\setminus\{0\}$ of the occupation measure of the strategy $\varphi^n$ for the initial state $x$. Now
\begin{eqnarray*}
\EE_{x}^{\varphi^{n}}[T_0]&=&\sum_{y=1}^\infty \eta_x^{\varphi^{n}}(y)=\sum_{y=x}^{x+n}
\left(\frac{1}{2}\right)^{y-x}
+2^{x}=\frac{1-\left(\frac{1}{2}\right)^{n+1}}{\frac{1}{2}}+2^x
=2+2^x-\frac{1}{2^n},
\end{eqnarray*}
and hence,
\begin{eqnarray*}
w^\ast(x)=\sup_{\pi\in\Delta^{\rm All}}\EE_x^\pi[T_0]\ge \sup_{n\ge 0} \EE_{x}^{\varphi^{n}}[T_0] =\sup_{n\ge 0}\left\{2+2^x-\frac{1}{2^n}\right\}=2+2^x,
\end{eqnarray*}
where the first equality is by the definition of the function $w^\ast(\cdot)$, see (\ref{PiunovskiyZhang2023cEqn06}).

On the other hand, it was shown in the proof of (b) of Theorem \ref{PiunovskiyZhang2023cThm01} for each $x\ge 1$ that $w^\ast(x)\le 2+2^x$, see (\ref{PiunovskiyZhang2023cEqn07}). Combining the previous two inequalities yields $w^\ast(x)=2+2^x$ for all $x\ge 1$. Thus, the statement is proved. $\hfill\Box$

Now we are in position to prove Proposition \ref{PiunovskiyZhang2023cThm02}.

\par\noindent\textit{Proof of Proposition \ref{PiunovskiyZhang2023cThm02}.}
(a) Recall that
\begin{eqnarray*}
v^\ast(x):=\inf_{\pi\in \Delta^{\rm All}}v^\pi(x)
\end{eqnarray*}
for each $x\in {\bf X}$. By definition of the function $w^\ast(\cdot)$, see (\ref{PiunovskiyZhang2023cEqn06}), $v^\ast(x)=-w^\ast(x)$ for each $x\in {\bf X}$. Now, according to Lemma \ref{PiunovskiyZhang2023cLemma01},
\begin{eqnarray}\label{PiunovskiyZhang2023cEqn15}
&&v^\ast(0)=0,\nonumber\\
&&v^\ast(x)=-2-2^x,~\forall~x\in\{1,2,\dots\}.
\end{eqnarray}

Suppose for contradiction that $\pi=\{\pi_n\}_{n=1}^\infty$ is an optimal strategy for the initial state $1$. If $\pi_1(1|1)=\varepsilon$ and $\pi_1(2|1)=1-\varepsilon$, then
\begin{eqnarray*}
v^\pi(1)=-1+\varepsilon\cdot\frac{1}{2}\cdot v^{{}^{(1,1)}\pi}(1)+(1-\varepsilon)\cdot\frac{1}{2}\cdot  v^{{}^{(1,2)}\pi}(2).
\end{eqnarray*}
where for $j\in\{1,2\}$  ${}^{(1,j)}\pi$ is the shifted strategy defined in (\ref{PiunovskiyZhang2023cEqn21}).
From (\ref{PiunovskiyZhang2023cEqn15}), we then see
 \begin{eqnarray*}
 v^\pi(1)\ge -1+\varepsilon\cdot\frac{1}{2}\cdot v^*_1+(1-\varepsilon)\cdot\frac{1}{2}\cdot v^*_2=-4+\varepsilon.
 \end{eqnarray*}
Hence, in order to obtain $v^\pi(1)=v^\ast(1)=-4$, it is necessary that $\varepsilon=0$.  The similar calculation for $\pi_2(a|x_0,a_1,x_1)$ with $x_0=1$, $a_1=2$ and $x_1=2$ leads to the equality $\pi_2(2|1,2,2)=1$, and so on. Thus,
\begin{eqnarray*}
v^\pi(1)=v^\varphi(1),
\end{eqnarray*}
where $\varphi$ is the deterministic stationary strategy given by $\varphi(x)\equiv 2$.

On the other hand,
\begin{eqnarray*}
v^\varphi(1)=-1-\frac{1}{2}-\frac{1}{4}-\ldots=-2>-4=v^\ast(1),
\end{eqnarray*}
and thus $v^\pi(1)>v^\ast(1)$, which is a desired contradiction against that $\pi$ is optimal for the initial state $1$. Consequently, there are no optimal strategies to problem (\ref{PiunovskiyZhang2023cEqn11}) for the MDP model with the given initial state $1$ and the given cost function $c(\cdot)$.

(b) The last assertion of this part is trivially true. In particular, by Lemma \ref{PiunovskiyZhang2023cLemma01},
\begin{eqnarray*}
\sup_{\pi\in\Delta^{\rm All}}\left\{\EE^\pi_{1}\left[\sum_{t=1}^\infty c^-(X_{t-1},A_t)\right]\right\}=w^\ast(1)=4<\infty,
\end{eqnarray*}
where the first equality is by (\ref{PiunovskiyZhang2023cEqn06}).
Thus, it follows from (a) and Proposition \ref{t0091} that condition (\ref{e1}) cannot be satisfied.

To show explicitly that Condition (\ref{e1}) is violated, consider again the deterministic stationary strategies $\varphi^n$ given by (\ref{PiunovskiyZhang2023cEqn08}), i.e.,
\begin{eqnarray*}
\varphi^n(x):= 2\cdot\I\{x\le n\}+\I\{x>n\}~\forall~ n=0,1,2,\ldots.
\end{eqnarray*}
Now
\begin{eqnarray*}
&& \inf_{N\ge n}\inf_{\pi\in\Delta^{\rm All}}\sum_{t=n+1}^N \EE^\pi_{1}\left[c(X_{t-1},A_t)\right]\le \EE^{\varphi^n}_1\left[c(X_{n},A_{n+1})\right]\\
&= & \left(\frac{1}{2}\right)^n\left[-1-(1-p_{n+1})-(1-p_{n+1})^2-\ldots\right]=- \left(\frac{1}{2}\right)^n\frac{1}{p_{n+1}}=-2,
\end{eqnarray*}
where the inequality holds because $c(\cdot)\le 0$, and the equalities hold by the similar argument as in (\ref{e33}). $\hfill\Box$

\begin{remark}
A strategy $\pi$ is called uniformly optimal if $v^\pi(x)=v^\ast(x)$ for all $x\in{\bf X}$. Consider an MDP model with a cost function $c(\cdot)$ such that $\sup_{\pi\in\Delta^{\rm All}} \EE^\pi_x\left[\sum_{n=0}^\infty c^-(X_n,A_{n+1})\right]<\infty$ and $v^\ast(x)\in (-\infty,\infty)$ for all $x\in {\bf X}$, where $v^\ast(\cdot)$ is defined by (\ref{PiunovskiyZhang2023cEqn17}). Then according to \cite[Theorem 2.2]{YuChi}, a deterministic stationary strategy $\varphi$ is uniformly optimal if and only if
\begin{equation}\label{e4}
v^*(x)= c(x,\varphi(x))+\int_{\bf X} v^*(y) p(dy|x,\varphi(x));~
\lim_{n\to\infty} \EE^\varphi_x\left[ v^*(X_n)\right]= 0.
\end{equation}
The two equalities in (\ref{e4}) are called the Dubins-Savage conditions.

Now, consider the optimal control problem (\ref{PiunovskiyZhang2023cEqn11}) for the MDP model in Example \ref{PiunovskiyZhang2023cExample01} with the cost function $c(x,a)\equiv -1$ for $x\in\{1,2,\dots\}$ and $c(0,a)\equiv 0.$ We have seen from Proposition \ref{PiunovskiyZhang2023cThm02}(a) that there is no uniformly optimal strategy. In particular, $\varphi(x)\equiv 2$ is not uniformly optimal. Let us verify this fact again by checking the Dubins-Savage conditions. They are sufficient and necessary for the uniform optimality of $\varphi$ because $c(\cdot)$ is $(-\infty,0]$-valued and $v^\ast(\cdot)$ is finite-valued. Now, we observe that $\varphi$ actually satisfies the first equality in (\ref{e4}). Nevertheless, the second equality in (\ref{e4}) is violated because
\begin{eqnarray*}
\EE_1^\varphi[v^*(X_n)]=\left(\frac{1}{2}\right)^n[-2-2^{n+1}]=-\left(\frac{1}{2}\right)^{n-1}-2\to -2~\mbox{ as } n\to\infty.
\end{eqnarray*}
Here $\left(\frac{1}{2}\right)^n$ is the probability that $X_n=n+1$.
\end{remark}

\section{Conclusion}\label{PiunovskiyZhang2023cSectionConclusion}
In conclusion, we showed that for a semi-continuous absorbing MDP with a fixed initial distribution, the continuity of the projection mapping $O$ from the space of strategic measures ${\cal P}$ to the space of occupation measures ${\cal D}$, both of which are endowed with their weak topologies, is equivalent to the MDP model being uniformly absorbing. This is confirmed  by an example.  Provided that the absorbing MDP model is semi-continuous,  the continuity of $O$ is a sufficient condition for the compactness of ${\cal D}$ in the weak topology. Whether it is also necessary is an interesting open problem for future studies.

Finally, we mention that if the MDP model is not absorbing, then ${\cal D}$ contains measures that are not (totally) finite. In this case, a different topology was introduced on ${\cal D}$ in e.g., \cite{DufourHoriguchiPiunovskiy:2012,PiunovskiyZhang:2023b}. In that topology, the projection mapping $O$ is continuous by definition.

\section*{Statements and Declarations} The authors declare that no funds, grants, or other support were received during the preparation of this manuscript. The authors have no relevant financial or non-financial interests to disclose. Both authors read and approved the final manuscript.   


\begin{thebibliography}{999}
\bibitem{b13} Altman, E. (1999). {\em Constrained Markov Decision Processes}. Chapman and Hall/CRC, Boca Raton.

\bibitem{b7} Bertsekas, D. and Shreve, S. (1978). {\em Stochastic Optimal Control}. Academic Press, New York.



\bibitem{Bogachev:2007} Bogachev, V. (2007). {\em Measure Theory (Volume 2)}. Springer, Berlin.

\bibitem{Bogachev:2018} Bogachev, V. (2018). {\em Weak Convergence of Measures}. American Mathematical Society, Providence.

\bibitem{Cavazos:1992} Cavazos-Cadena, R. and Hern\'{a}ndez-Lerma. (1992). Equivalence of Lyapunov stability criteria in a
class of Markov decision processes. {\em Appl. Math. Optim.} {\bfseries 26}, 113--137.

\bibitem{Derman:1970} Derman, C. (1970). {\em Finite State Markovian Decision Processes}. Academic Press, New York.



\bibitem{Dufour:2020} Dufour, F. and Genadot, A. (2020). On the expected total cost with unbounded returns for Markov decision processes. {\em Appl. Math. Optim.} {\bfseries 82}, 433--450.

\bibitem{DufourHoriguchiPiunovskiy:2012} Dufour, F., Horiguchi, M. and Piunovskiy, A. (2012). The expected total cost criterion for Markov decision processes under constraints: a convex analytic approach. {\em Adv. Appl. Prob.} {\bfseries 44}, 774--793.



\bibitem{b114} Feinberg, E.A. and Rothblum, U. (2012). Splitting randomized stationary policies in total-reward Markov decision processes. {\em Math. Oper. Res.} \textbf{37}, 129--153.

\bibitem{b16} Feinberg, E.A. and Piunovskiy, A. (2019). Sufficiency of deterministic policies for atomless discounted and uniformly absorbing MDPs with multiple criteria. {\em SIAM J. Control Optim.}  \textbf{57}, 163--191.


\bibitem{b12} Hern\'{a}ndez-Lerma, O. and Lasserre, J.B. (1999). \emph{Further Topics on Discrete-Time Markov Control Processes}. Springer-Verlag, New York.


\bibitem{Hordijk:1975} Hordijk, A. (1975). {\em Dynamic Programming and Potential Theory}. Mathematish Centrum, Amsterdam.


\bibitem{Collins:2006} James, H. and Collins, E. (2006). An analysis of transient Markov decision processes. {\em J. Appl. Probab.} {\bfseries 43}, 603--621.

\bibitem{Kurushima:2018} Kurushima, A., Piunovskiy, A. and Zhang, Y. (2018). Nowak's theorem on probability measures induced by strategies revisited. {\em Theory Probab. its Appl.} {\bfseries 62}, 328--334.

\bibitem{Nowak:1988} Nowak, A. (1988). On the weak topology on a space of probability measures induced by policies. {\em Bulletin Polish Acad. Sci. Math.} {\bfseries 36}, 181--186.


\bibitem{PiunovskiyZhang:2023b} Piunovskiy, A. and Zhang, Y. (2023).  Extreme occupation measures in Markov decision processes with an absorbing state. {\em SIAM J. Control Optim.}, in press.


\bibitem{b18} Sch\"al, M. (1975). On dynamic programming: compactness of the space of policies. {\em Stoch. Proc. Appl.} {\bfseries 3}, 345--364.

\bibitem{YuChi}  Yushkevich, A.A. and Chitashvili, R. Ya. (1982). Controllable random sequences and Markov chains. {\em Russian Math. Surveys} {\bfseries 37}, 239-–274.


\end{thebibliography}
\end{document}